\documentstyle[12pt]{article}
\newtheorem{t.}{Theorem}[section]
\newtheorem{d.}[t.]{Definition}
\newtheorem{l.}[t.]{Lemma}
\newtheorem{r.}[t.]{Remark}
\newtheorem{c.}[t.]{Corollary}
\newtheorem{e.}[t.]{Example}
\begin{document}
\title {On (Co)Homology of Triangular Banach Algebras \footnote{
This research was in part supported by a grant from IPM (No. 82460027).}}
\author{{\bf Mohammad Sal Moslehian} \\ Dept. of Math., Ferdowsi Univ. P.O.Box 1159, Mashhad 91775\\and\\Institute for Studies in Theoretical Physics and Mathematics (IPM)\\ Iran\\ E-mail: msalm@math.um.ac.ir\\Home: http://www.um.ac.ir/$\sim$moslehian/}
\date{}
\maketitle
\font\tt=cmtt10
\begin{abstract}
Suppose that $A$ and $B$ are unital Banach algebras with units $1_A$ and $1_B$, respectively, $M$ is a unital Banach $A-B$-bimodule, ${\cal T}= \left [\begin{array}{cc}A&M\\0&B \end{array}\right]$ is the triangular Banach algebra, $X$ is a unital ${\cal T}$-bimodule, $X_{AA}=1_AX1_A$, $X_{BB}=1_BX1_B$, $X_{AB}=1_AX1_B$ and $X_{BA}=1_BX1_A$. Applying two nice long exact sequences related to $A$, $B$, ${\cal T}$, $X$, $X_{AA}$, $X_{BB}$, $X_{AB}$ and $X_{BA}$ we establish some results on (co)homology of triangular Banach algebras.
\end{abstract}
\newpage

\section{Introduction}

Topological homology arose from the problems concerning extensions by H. Kamowitz who introduced the Banach version of Hochschild cohomology groups in 1962 [11], derivations by R. V. Kadison and J. R. Ringrose [9, 10] and amenability by B.E. Johnson [8] and has been extensively developed by A. Ya. Helemskii and his school. In addition, this area includes a lot of problems concerning automorphisms, fixed point theorems, perturbations, invariant means, topology of spectrum, ... . [6].

This article deals with the cohomology and homology of triangular Banach algebras, i.e. algebras of the form ${\cal T}= \left [\begin{array}{cc}A&M\\0&B \end{array}\right]$ in which $A$ and $B$ are unital Banach algebras and $M$ is a unital Banach $A-B$-bimodule. These algebras introduced by Forrest and Marcoux [1], motivated by work of Gilfeather and Smith in [4]. Forrest and Marcoux also studied and directly computed some cohomology groups of triangular Banach algerbras (see [2] and [3]).\\
In this paper, after preliminaries as a foundation for our work, we present two long exact sequences and apply them to give some significant isomorphisms and vanishing theorems.

\section{Preliminaries}

We begin with some observations concerning cohomology and homology of Banach algebras. Some sources of references to the subject are [6] and [7].

Let {\bf Lin} denote the category of linear spaces and linear operators. A sequence $\cdots \longleftarrow X_n \stackrel{d_n}{\longleftarrow}X_{n+1} \longleftarrow \cdots ~~,~~ {\cal X} = \{X,d\}$ (resp. $\cdots \longrightarrow X^n \stackrel{\delta^n}{\longrightarrow}X^{n+1} \longrightarrow \cdots~~,~~ {\cal X} = \{X,\delta\}$) in a subcategory of ${\bf Lin}$ is said to be a (chain) complex (resp. (cochain) complex) if $d_{n-1} \circ d_n = 0$ (resp. $\delta^n \circ \delta^{n-1} = 0$).

Suppose that $A$ is a Banach algebra and $X$ is a Banach $A$-bimodule.

For $n = 0,1,2, \cdots $, let $ C^n(A,X)$ be the Banach space of all bounded n-linear mappings from $A \times \cdots \times A$ into $X$ together with multilinear operator norm $\parallel f \parallel = \sup \{ \parallel f(a_1, \ldots ,a_n)\parallel ;~ a_i \in A, \parallel a_i \parallel \leq 1,  1\leq  i\leq n \}$, and $C^0(A,X) = X$. The elements of $C^n(A,X)$ are called n-dimensional cochains. Consider the sequence\\ $0 \longrightarrow C^0(A,X)\stackrel{\delta^0}{\longrightarrow} C^1(A,X)\stackrel{\delta^1}{\longrightarrow}\cdots~~({\widetilde{C}(A,X)}$, where $\delta^0x(a) = ax - xa$ and for $~n = 0,1,2, \cdots, \delta^nf(a_1, \ldots, a_{n + 1}) = a_1f(a_2, \ldots, a_{n + 1}) \\+ \sum_ {k = 1}^n (-1)^kf(a_1, \ldots, a_{k - 1}, a_ka_{k + 1}, \ldots, a_{n + 1}) + (-1)^{n + 1}f(a_1, \ldots, a_n)\\a_{n + 1}~~( x \in X, ~a, a_1,
\ldots, a_{n + 1} \in A, ~f \in C^n(A,X))$.

It is straightforward to verify that the above sequence is a complex. $\widetilde{C}(A,X)$ is called the standard cohomology complex or Hochschild-Kamowitz complex for $A$ and $X$. The n-th cohomology group of $\widetilde{C}(A,X)$ is said to be n-dimensional (ordinary or Hochschild) cohomology group of $A$ with coefficients in $X$ and denoted by $H^n(A,X)$. The spaces $Ker\delta ^n$ and $Im\delta^{n-1}$ are denoted by $Z^n(A,X)$ and $B^n(A,X)$, and their elements are called $n$-dimensional cocycles and $n$-dimensional coboundaries, respectively. Hence $H^n(A,X) = Z^n(A,X)/B^n(A,X)$. Note that $H^n(A,X)$, generally speaking, is a complete seminormed space.

Assume that $C_0(A,X)= X$ and for $n = 1,2, \cdots, C_n(A,X) = {\underbrace{A\stackrel{\wedge}{\otimes}\cdots\stackrel{\wedge}{\otimes}A}_n} \stackrel{\wedge}{\otimes}X$ in which $\stackrel{\wedge}{\otimes}$ denotes the projective tensor product of Banach spaces. The elements of $C_n(A,X)$ are called $n$-dimensional chains. Consider the complex\\
$0 \longleftarrow C_0(A,X)\stackrel{d_0}{\longleftarrow} C_1(A,X)\stackrel{d_1}{\longleftarrow}\cdots~~(\widehat{C}(A,X))$;~ where $d_n(a_1\otimes \ldots \otimes a_{n+1}\otimes x) = a_2\otimes \ldots \otimes a_{n+1}\otimes a_1x + \sum_{k = 1}^n (-1)^k a_1\otimes \ldots \otimes a_ka_{k+1}\otimes \ldots \otimes a_{n+1}\otimes x + (-1)^{n+1}a_1\otimes \ldots \otimes a_n\otimes xa_{n+1}$.

The $n$-th homology group of $\widehat C (A,X)$ is called (ordinary) homology group of $A$ with coefficients in $X$. It is denoted by $H_n(A,X)$ which is a complete seminormed space. 

The dual $X^*$ of the Banach $A$-bimodule $X$ is again a Banach $A$-bimodule with respect to the following actions:
$$(af)(x)=f(xa) ~~,~~ (fa)(x)=f(ax)~~;~~ f\in X^* , a\in A, x\in X$$
In particular, $A^*$ is a Banach bimodule over $A$. 

A complex ${\cal X}=\{X,d\}$ in a category of Banach modules is called admissible if it splits as a complex of Banach spaces and continuous linear operators, i.e. the kernels of all its morphisms are topologically complemented.

An additive functor $F$ is said to be exact if for every admissible complex ${\cal X}=\{X,d\}$ the complex $\underline{F}({\cal X})=\{F(X),F(d)\}$ is exact in the category {\bf Lin}. Notice that $\underline{F}$ is a functor.

A unital left Banach module $P$ over a unital Banach algebra $A$ is said to be projective if the functor $_Ah(P,?)$ is exact. Recall that for left Banach $A$-modules $X$ and $Y,~ _Ah(?,?)$ takes left $A$-modules $X$ and $Y$ to $_Ah(X,Y)=\{f:X\to Y; f$ {\rm is a bounded left} A-{\rm module map} $\}$. Indeed $_Ah(?,?)$ is a bifunctor contravariant in the first variable and covariant in the second.

A left Banach $A$-module (resp. right Banach $A$-module, Banach $A-B$-bimodule) $X$ is called projective if $X$ as a left Banach unital $A_+$-module (resp. left Banach unital $A_+^{op}$-module, left Banach unital $A_+\stackrel{\wedge}{\otimes}B_+^{op}$-module) is projective, where $A_+=A\oplus {\bf C}$ denotes the unitization of the Banach algebra $A. A^{op}$, the so-called opposite to $A$, is the space $A$ equipped with the multiplication $a\circ b=ba$.

A complex $0 \gets X_0 \stackrel{d_0}{\gets}X_1 \gets \cdots ~~,~~ ({\cal X})$ is called a resolution of the $A$-module $X$ if the complex $0 \gets X \stackrel{\varepsilon}{\gets} X_0 \stackrel{d_0}{\gets}X_1 \gets \cdots~~$ is admissible. By a projective resolution we mean one in which the $X_i$'s are projective. 

Every left $A$-module $X$ admits sufficiently many projective resolutions, especially it admits the normalized bar-resolution ${\cal B}(X)$ as follows:

Consider free modules $B_n(X)=A_+  \stackrel{\wedge}{\otimes}(\underbrace{A \stackrel{\wedge}{\otimes} \cdots  \stackrel{\wedge}{\otimes}A}_n \stackrel{\wedge}{\otimes}X)$ and $B_0(X)=A_+ \stackrel{\wedge}{\otimes}X$, and operators $\pi : A_+\stackrel{\wedge}{\otimes}X \to X$ and $d_n : B_{n + 1}(X) \to B_n(X)$, well-defined by $\pi(a \otimes x)=ax,~~
d_n(a\otimes a_1\otimes \ldots \otimes a_{n+1}\otimes x)=aa_1\otimes a_2\otimes \ldots \otimes a_{n+1}\otimes x + \sum_{k=1}^n (-1)^k a\otimes a_1\otimes \ldots \otimes a_ka_{k+1}\otimes \ldots \otimes a_{n+1}\otimes x + (-1)^{n+1}a\otimes a_1\otimes \ldots \otimes a_n\otimes a_{n+1}x$.\\
Then the sequence $0 \gets X \stackrel{\pi}{\gets}B_0(X) \stackrel{d_0}{\gets}B_1(X)\stackrel{d_1}{\gets} \cdots$ is a projective resolution of $X$.\\
We denote the complex $0 \gets B_0(X) \stackrel{d_0}{\gets}B_1(X)\stackrel{d_1}{\gets} \cdots$ by ${\cal B}(X)$. In fact ${\cal B}$ induces a functor.

Let $F$ be an additive functor. Then the functor $F_n=H_n \circ \underline{F} \circ {\cal B}$ is called the $n$-th projective derived functor of $F. ~F_n$ is independent of the choice of resulotion. The projective derived cofunctors could be defined in a similar way. 

For a left Banach $A$-module $Y$, let $Ext_A^n(?,Y)$ denote the $n$-th projective derived cofunctor of $_Ah(?,Y)$. Given a right Banach $A$-module $X$, denote the $n$-th projective derived functor of $X \stackrel{\wedge}{\otimes}_A?$ by $Tor_n^A(X,?)$. Recall that for a right Banach $A$-module $X$ and a left Banach $A$-module $Y$ the projective tensor product of modules $X$ and $Y$ is defined to be the quotient space $X \stackrel{\wedge}{\otimes}_AY=(X \stackrel{\wedge}{\otimes}Y)/L$ where $L$ denotes the closed linear span of all elements of the form $xa \otimes y - x \otimes ay~ ( x \in X , a \in A , y \in Y )$. In fact $? \stackrel{\wedge}{\otimes}_A?$ is a bifunctor covariant in both variables.

\section{Main results}

Suppose that $A$ and $B$ are unital Banach algebras with units $1_A$ and $1_B$, and Banach space $M$ is a unital Banach $A-B$-bimodule. Then ${\cal T}= \left [\begin{array}{cc}A&M\\0&B \end{array}\right]=\{ \left [\begin{array}{cc}a&m\\0&b \end{array}\right]; ~a\in A, m\in M, b\in B\}$ with the usual $2\times 2$ matrix addition and formal multiplication equipped with the norm $\| \left [\begin{array}{cc}a&m\\0&b \end{array}\right]\|=\| a\|+\| m\|+\| b\|$ is a Banach algebra which is called a triangular Banach algebra.[1]

Let $X$ is a unital Banach ${\cal T}$-bimodule, $X_{AA}=1_AX1_A, X_{BB}=1_BX1_B, X_{AB}=1_AX1_B$ and $X_{BA}=1_BX1_A$. 

Applying homological techniques we could establish the following long exact sequences (see [5]):\\

$0 \stackrel{\pi^{-1}}{\to} H^0({\cal T},X) \stackrel{\phi^0}{\to} H^0(A,X_{AA})\oplus H^0(B,X_{BB}) \stackrel{\delta^0}{\to} Ext^0_{A\stackrel{\wedge}{\otimes} B^{op}}(M,X_{AB})\\
\stackrel{\pi^0}{\to} H^1({\cal T},X) \stackrel{\phi^1}{\to} H^1(A,X_{AA})\oplus H^1(B,X_{BB}) \stackrel{\delta^1}{\to} Ext^1_{A\stackrel{\wedge}{\otimes} B^{op}}(M,X_{AB})\\
\stackrel{\pi^1}{\to} H^2({\cal T},X) \stackrel{\phi^2}{\to} H^2(A,X_{AA})\oplus H^2(B,X_{BB}) \stackrel{\delta^2}{\to} Ext^2_{A\stackrel{\wedge}{\otimes} B^{op}}(M,X_{AB}) \to \cdots$
\\
~\\
$\cdots \stackrel{\rho_2}{\to} Tor_2^{A\stackrel{\wedge}{\otimes} B^{op}}(M,X_{BA})\stackrel{d_1}{\to} H_2(A,X_{AA})\oplus H_2(B,X_{BB}) \stackrel{\psi_2}{\to}H_2({\cal T},X)\\
\stackrel{\rho_1}{\to} Tor_1^{A\stackrel{\wedge}{\otimes} B^{op}}(M,X_{BA})\stackrel{d_1}{\to} H_1(A,X_{AA})\oplus H_1(B,X_{BB}) \stackrel{\psi_1}{\to}H_1({\cal T},X)\\
\stackrel{\rho_0}\to Tor_0^{A\stackrel{\wedge}{\otimes} B^{op}}(M,X_{BA})\stackrel{d_0}{\to}
H_0(A,X_{AA})\oplus H_0(B,X_{BB}) \stackrel{\psi_0}{\to} H_0({\cal T},X)\stackrel{\rho_{-1}}{\to}0$
\\

Using these nice sequences we shall obtain some significant results:

\begin{t.}
Let $X_{AB}=0$ and ${\cal T}= \left [\begin{array}{cc}A&M\\0&B \end{array}\right]$. Then $H^n({\cal T}, X)\simeq H^n(A,X_{AA})\oplus H^n(B,X_{BB})$ for all $n\geq 0$.
\end{t.}
{\bf Proof.} If $X_{AB}=0$, then $Ext^{n-1}_{A\stackrel{\wedge}{\otimes} B^{op}}(M,X_{AB})=Ext^n_{A\stackrel{\wedge}{\otimes} B^{op}}(M,X_{AB})=0$.\\
Hence $0 \stackrel{\pi^{n-1}}{\to} H^n({\cal T},X) \stackrel{\phi^n}{\to} H^n(A,X_{AA})\oplus H^n(B,X_{BB}) \stackrel{\delta^n}{\to} 0.\Box$

\begin{c.}
${\cal T}= \left [\begin{array}{cc}A&M\\0&B \end{array}\right]$ is weakly amenable iff so are $A$ and $B$.{\rm [2, Corollary 3.5]}
\end{c.}
{\bf Proof.} $X={\cal T^*}$ is a Banach ${\cal T}$-bimodule for which clearly $X_{AA}=A^*, X_{BB}=B^*, X_{AB}=0$ and $X_{BA}=M^*$. Then the previous theorem, with $n=1$, implies that
$$H^1({\cal T},{\cal T}^*)=H^1(A,A^*)\oplus H^1(B,B^*).$$
Hence ${\cal T}$ is weakly amenable iff so are $A$ and $B.\Box$

\begin{t.}
Let $X_{BA}=0$ and ${\cal T}= \left [\begin{array}{cc}A&M\\0&B \end{array}\right]$. Then $H_n({\cal T}, X)\simeq H_n(A,X_{AA})\oplus H_n(B,X_{BB})$ for all $n\geq 0$.
\end{t.}
{\bf Proof.} If $X_{BA}=0$, then $Tor_{n-1}^{A\stackrel{\wedge}{\otimes} B^{op}}(M,X_{AB})=Tor_n^{A\stackrel{\wedge}{\otimes} B^{op}}(M,X_{AB})=0$.\\
Hence $0 \stackrel{d_n}{\to} H_n(A,X_{AA})\oplus H_n(B,X_{BB})\stackrel{\psi_n}{\to}H_n({\cal T},X)\stackrel{\rho_{n-1}}{\to} 0.\Box$ 

\begin{c.}
$H_n({\cal T}, {\cal T})\simeq H_n(A,A)\oplus H_n(B,B)$ if ${\cal T}= \left [\begin{array}{cc}A&M\\0&B \end{array}\right]$.
\end{c.}
{\bf Proof.} For $X={\cal T}$ we have $X_{AA}=A, X_{BB}=B, X_{AB}=M$ and $X_{BA}=0.\Box$

\begin{c.}
If ${\cal T}= \left [\begin{array}{cc}A&M\\0&B \end{array}\right]$, then $H_n({\cal T}, M)=0$. In particular, with ${\cal T}_m= \left [\begin{array}{cc}A&{\cal T}_{m-1}\\0&B \end{array}\right]$ and ${\cal T}_0={\cal T}$, we conclude that $H_n({\cal T}, {\cal T}_m)=0$.
\end{c.}
{\bf Proof.} Note that $X_{AA}=A, X_{BB}=B, X_{AB}=M$ and $X_{BA}=0$, if $X=M.\Box$

\begin{t.}
Denote by $\tau(D)$ the set of all bounded traces over the Banach algebra $D$, i.e. $\tau(D)= \{f\in D^*; f(d_1d_2)=f(d_2d_1), {\rm ~for ~all~} d_1,d_2\in D\}$. Then 
$\tau({\cal T}) \simeq \tau(A)\oplus \tau(B)$, if ${\cal T}= \left [\begin{array}{cc}A&M\\0&B \end{array}\right]$.
\end{t.}

{\bf Proof.} $0 \to H^0({\cal T},{\cal T}^*) \stackrel{\phi^0}{\to} H^0(A,A^*)\oplus H^0(B,B^*) \stackrel{\delta^0}{\to}\\ Ext^0_{A\stackrel{\wedge}{\otimes} B^{op}}(M,{\cal T}^*_{AB})=0$, since ${\cal T}^*_{AB}=0$. Note then that $H^0(D,D^*)=\tau(D).\Box$

\begin{r.} {\rm Thanks to Niels Jakob Laustsen for his comment on the fact that there is a direct proof for Theorem 3.6:\\
The equality $\left [\begin{array}{cc}0&m\\0&0 \end{array}\right]=\left [\begin{array}{cc}1&0\\0&0 \end{array}\right]\left [\begin{array}{cc}0&m\\0&0 \end{array}\right]-\left [\begin{array}{cc}0&m\\0&0 \end{array}\right]\left [\begin{array}{cc}1&0\\0&0 \end{array}\right]$ implies that $f(\left [\begin{array}{cc}0&m\\0&0 \end{array}\right]=0$ for every $f\in \tau({\cal T})$. Then $f_A(a)=f(\left [\begin{array}{cc}a&0\\0&0 \end{array}\right])$ and $f_B(b)=f(\left [\begin{array}{cc}0&0\\0&b \end{array}\right])$ give two bounded traces over $A$ and $B$ respectively. Conversely, if we have two bounded traces $f_1$ and $f_2$ on $A$ and $B$, resp., then $f(\left [\begin{array}{cc}a&m\\0&b \end{array}\right]=f_1(a)+f_2(b)$ defines a bounded trace over ${\cal T}$.}\end{r.}

\begin{t.}
Let $A$ be a unital Banach algebra with $H^n(A,A)=0$ for some $n>1$, and $M$ be a left Banach $A$-module, then $H^n({\cal T},{\cal T})\simeq H^n(A,B(M))$ in which ${\cal T}=\left [\begin{array}{cc}A&M\\0&{\bf C} \end{array}\right ]$. {\rm ( See [3, Corollary 4.2])}\end{t.}

{\bf Proof.} Put $X={\cal T}$. By $Ext^n_A(X,Y)\simeq H^n(A,B(X,Y))$ and $H^n({\bf C},{\bf C})=0$, the exact sequence 

$\cdots \stackrel{\phi^{n-1}}{\to} H^{n-1}(A,A)\oplus H^{n-1}({\bf C},{\bf C}) \stackrel{\delta^{n-1}}{\to} Ext^{n-1}_{A\otimes {\bf C}^{op}}(M,M) \\
\stackrel{\pi^{n-1}}{\to} H^n({\cal T},{\cal T}) \stackrel{\phi^n}{\to} H^n(A,A)\oplus H^n({\bf C},{\bf C}) \stackrel{\delta^n}{\to} \cdots$ \\
gives rise to $$\cdots \to 0 \to H^{n-1}(A,M) \to H^n({\cal T},{\cal T}) \to 0 \to \cdots.$$
Hence $H^n({\cal T},{\cal T})\simeq H^n(A,B(M)).\Box$

\begin{e.} {\rm Suppose that $A$ is a hyperfinite von Neumann algebra acting on a Hilbert space $H$. $M=H$ is a left $A$-module via $a.\xi=a(\xi), a\in A, \xi\in H$. It follows from [12, Corollary 3.4.6], $H^n(A,A)=0$. So $H^n(\left [\begin{array}{cc}A&H\\0&{\bf C} \end{array}\right ],\left [\begin{array}{cc}A&H\\0&{\bf C} \end{array}\right ])=H^{n-1}(A,B(H))$. In particular, $H^2(\left [\begin{array}{cc}A&H\\0&{\bf C} \end{array}\right ],\left [\begin{array}{cc}A&H\\0&{\bf C} \end{array}\right ])=H^1(A,B(H))=0$, by [12, Theorem 2.4.3]. (See [3, Example 4.2])\\}\end{e.}

\end{document}